\documentclass{amsart}

\usepackage[foot]{amsaddr}
\usepackage{tikz}
\usepackage{natbib}
\bibpunct{(}{)}{;}{a}{}{,}
\renewcommand{\cite}{\citep}

\title%[There's Glory for You!]{There's Glory for You!\\ 
{Redefining Revolutions}
\author[Andrew Aberdein]{Andrew Aberdein$^{*}$}
\address{$^{*}$School of Arts \& Communication, Florida Institute of Technology, Melbourne FL.}
\thanks{Published as: Andrew Aberdein, Redefining revolutions, In \emph{The Kuhnian Image of Science: Time for a Decisive Transformation?} (Moti Mizrahi, ed.), Rowman and Littlefield, London, 2018, pp.~133--154.}
\begin{document}
\maketitle

%In their account of theory change in logic, Aberdein and Read distinguish `glorious' from `inglorious' revolutions---only the former preserves all `the key components of a theory' \cite{Aberdein09a}.
%A widespread view, expressed in these terms, is that empirical science characteristically exhibits inglorious revolutions but that revolutions in mathematics are at most glorious \cite{Gillies92}.
%Here are three possible responses:
%\begin{enumerate}\addtocounter{enumi}{-1}
% 	\item Accept that empirical science and mathematics are methodologically discontinuous;
%	\item Argue that mathematics can exhibit inglorious revolutions;
%	\item Deny that inglorious revolutions are characteristic of science.
%\end{enumerate}
%Where Aberdein and Read take option 1, option 2 is preferred by Mizrahi \cite{Mizrahi14b}.
%This paper seeks to resolve this disagreement through consideration of some putative mathematical revolutions.

\section{A Nice Knockdown Argument}\label{Sec:Mizrahi}%Mizrahi against incommensurability}
Moti Mizrahi does an admirable job in pruning the thicket that has grown up around Thomas Kuhn's incommensurability thesis. He begins with a distinction between two versions of the thesis:
\begin{quotation}
\begin{description}
\item[Taxonomic Incommensurability (TI)] Periods of scientific change (in particular, revolutionary change) that exhibit TI are scientific developments in which existing concepts are replaced with new concepts that are incompatible with the older concepts. The new concepts are incompatible with the old concepts in the following sense: two competing scientific theories are conceptually incompatible (or incommensurable) just in case they do not share the same ``lexical taxonomy''. A lexical taxonomy contains the structures and vocabulary that are used to state a theory. %(See Kuhn 2000, pp. 14Ð5, 63, 92Ð7, 229Ð33, 238Ð9, and 242Ð4; cf. Sankey 1997)2
\item[Methodological Incommensurability (MI)] There are no objective criteria of theory evaluation. The familiar criteria of evaluation, such as simplicity and fruitfulness, are not a fixed set of rules. Rather, they vary with the currently dominant paradigm. %(See Kuhn 1962, p. 94, 1970, p. 200, 1977, pp. 322, 331; cf. Sankey and Hoyningen-Huene 2001, p. xiii)3
\cite[362; references omitted]{Mizrahi14b}.
\end{description}
\end{quotation}
Mizrahi's focus is exclusively on (TI), a focus that I will share. He proceeds to argue that, understood as (TI), the incommensurability thesis is poorly motivated. His critique differs from that of many other authors in focussing on the weakness of the arguments in support of (TI), rather than the strength of the arguments against it. He claims that the former arguments must be either deductive or inductive. In both cases, he presents a counterargument. Against deductive support, he argues as follows:
\begin{enumerate}
\item Reference change (discontinuity) is conclusive evidence for (TI) only if reference change (discontinuity) entails incompatibility of conceptual content.
\item Reference change (discontinuity) does not entail incompatibility of conceptual content. 

Therefore:
\item It is not the case that reference change (discontinuity) is conclusive evidence for (TI).
\cite[367]{Mizrahi14b}.
\end{enumerate}
Against inductive support, he argues as follows:
\begin{enumerate}
\item There is a strong inductive argument for (TI) only if there are no rebutting defeaters against (TI).
\item There are rebutting defeaters against (TI). 

Therefore:
\item It is not the case that there is a strong inductive argument for (TI).
\cite[371]{Mizrahi14b}.
\end{enumerate}

The merits of this critique have been addressed elsewhere \cite{Kindi15,Marcum15,Patton15,Mizrahi15e,Mizrahi15d}.
In this chapter, I will take a somewhat different tack, by examining the implications for the philosophy of mathematical practice, specifically the debate whether there can be mathematical revolutions. Hence I will not engage closely with all the details of Mizrahi's argument. However, I do wish to draw attention to his use of `rebutting defeater'. Before introducing an example from the history of medicine of a revolution in which conceptual continuity is displayed, he stresses that
\begin{quotation}
the following episode is not supposed to be a counterexample against (TI). It is not meant to be a refutation of (TI). Rather, it shows that an inductive argument based on a few selected historical episodes of scientific change does not provide strong inductive support for (TI). Or, to put it another way, this episode---and others like it---counts as what Pollock %(1987, pp. 481Ð518) 
calls a \emph{rebutting defeater}, i.e.\ a \emph{prima facie} reason to believe the negation of the original conclusion; in this case, the negation of (TI) \cite[368; reference omitted]{Mizrahi14b}.
\end{quotation}
That is, cases of revolutionary change without conceptual discontinuity are rebutting defeaters for (TI) since they are reasons to think that we can have one without the other. They are not undercutting defeaters, since the familiar cases of revolutionary change \emph{with} conceptual discontinuity are still reasons to believe (TI), but if we have as many reasons to disbelieve it as to believe it, we should probably suspend our judgment.

%\cite{DAgostino14}.
%\cite{Balaguer09}.
%\cite{Bledin08}.
%\cite{Rusnock12}.
%\cite{Swiggers82}.
%\cite{Collier84}.
%\cite{Faye13}.
%`conceptual change in empirical theories \dots\ Only the most extreme type---replacement of dimensions---comes close to a revolution' \cite[1039]{Gardenfors13}.

\section{There's Glory for You!}
Alice's encounter with Humpty Dumpty is well-known to philosophers:
\begin{quotation}
``I don't know what you mean by `glory,'\,'' Alice said. 

Humpty Dumpty smiled contemptuously. ``Of course you don't---till I tell you. I meant `there's a nice knock-down argument for you!'\,'' 

``But `glory' doesn't mean `a nice knock-down argument,'\,'' Alice objected. 

``When \emph{I} use a word,'' Humpty Dumpty said, in rather a scornful tone, ``it means just what I choose it to mean---neither more nor less.'' 

``The question is,'' said Alice, ``whether you \emph{can} make words mean so many different things.'' 

``The question is,'' said Humpty Dumpty, ``which is to be master---that's all.'' \cite[106 f.]{Carroll97}.
\end{quotation}
%\cite{Ketland12}.
%\cite{Schupbach15}.
Several echoes of this passage may be heard in this chapter, most clearly in a Humpty-ish passage of my own: 
\begin{quote}
A \textit{glorious} revolution occurs when the key components of a theory are preserved, despite changes in their character and relative significance. (We will refer to such preservation, constitutive of a glorious revolution, as \textit{glory}.) An \textit{inglorious} revolution occurs when some key component(s) are lost, and perhaps other novel material is introduced by way of replacement. \dots\ A \textit{paraglorious} revolution occurs when all the key components are preserved, as in a glorious revolution, but new key components are also added. 
\dots\ A \textit{null} revolution\footnote{Less happily, I also referred to theories in null revolution as being in stasis, but `stasis' is a false friend: although it has come to mean an absence of change, for Aristotle it meant something much like revolution, and is often translated as such \cite[18]{Howell85}.}
\dots\ %, as it were) 
when none of its key components change at all \cite[618 f.]{Aberdein09a}.
\end{quote}

%In the interests of clarity, h
Here is a slightly more formal characterization of this fourfold distinction. Let us specify that the key components of a theory $T_{n}$ comprise a structure of some sort, $K_{n}$. Then succession between theories $T_{n}$ and $T_{n+1}$
would be a null revolution if their respective key components are unchanged, that is $K_{n}= K_{n+1}$, and
a glorious revolution if the key components are in some sense isomorphic, $K_{n}\cong K_{n+1}$, that is if there is a well-motivated bijection between them which respects their roles in each theory.
By contrast, in a paraglorious revolution there would be an analogous isomorphic embedding of the key terms of the old theory within the new, $K_{n}\hookrightarrow K_{n+1}$, such that the new theory contains terms with no clear counterpart in the old.
And, in what we may call a strict inglorious revolution, there would be a similar isomorphic embedding of the key terms of the new theory within the old, $K_{n}\hookleftarrow K_{n+1}$, since there are components of the old theory which have been irretrievably lost in the new. %\footnote
{In other words, inglorious revolutions would exhibit `Kuhn loss', a loss of (actual or potential) explanatory power \cite[for further discussion, see][111 ff.]{Votsis11}.} %Heinz Post
(The more general sense of inglorious revolution may be thought of as a strict inglorious revolution combined with a paraglorious revolution. That is, there would be some structure $K_{n}^{\prime}$, not necessarily corresponding to any actually espoused theory, such that $K_{n}\hookleftarrow K_{n}^{\prime}\hookrightarrow K_{n+1}$.)

In my earlier presentation of this distinction, I addressed a number of questions, the most important of which are what components are `key' and how are they `preserved'?
The simplest answer to the first question would be to make \emph{all} components of a theory key, or at least all components without which the theory could not be articulated. A more subtle account would permit distinctions between the theory proper and auxiliary theories, tentative extensions, and other inessential components, but we need not explore that account here.%
\footnote{The sort of distinction I have in mind is that Stathis Psillos draws between idle and essentially contributing constituents \cite[S311]{Psillos96} or Philip Kitcher between presuppositional and working posits \cite[149]{Kitcher93}.}
Indeed, mathematical theories are less trouble than empirical theories in this respect: they characteristically have fewer components and their dependencies are much more clearly stated.
The second question is much more of a challenge. Indeed, it is central to understanding what makes a revolution revolutionary: how much change is required for a revolution? Conversely, how much change can a theory undergo without revolution? 
In other words, just what do we mean by `glory'?
An obvious starting point would be taxonomic commensurability, that is the absence of TI.
Notice, incidentally, that lexical taxonomy is shared across neither inglorious nor paraglorious revolutions. I have defined inglorious and paraglorious revolutions such that the distinction is essentially a matter of historic sequence: whether such a transition counts as inglorious or paraglorious will depend on which theory came first. The definition of (TI), however, despite references to `new' and `old', does not seem to directly appeal to chronology.
%I will return to this question in \S\ref{Sec:Con}.

Another question posed by this framework concerns the transitivity (or not) of these different characterizations of change. A succession of null revolutions must be a null revolution, because we have stipulated strict identity between key components. But a succession of glorious revolutions need not be a glorious revolution: a series of comparatively small changes might add up to a big change, as in a sorites sequence of small changes of colour from red to blue. 
Likewise, the preservation aspect of paraglorious revolution may fail over a long enough sequence of such revolutions, making the sequence as a whole inglorious.
Inglorious revolutions themselves are more straightforward: in principle, two consecutive inglorious revolutions might cancel each other out, so inglorious revolution must be intransitive.

A final concern leads us back to Humpty Dumpty: mere survival (or not) of vocabulary is not what is at issue.% 
\footnote{Mizrahi also addresses this issue \cite[374]{Mizrahi14b}.}
Conceptual change might be disguised by a shift in the meaning of a shared vocabulary; conversely, a drastic change of vocabulary may give a misleading impression of change when nothing substantive has occurred. This issue is familiar from political examples: Augustus strategically reused much of the terminology of the old Roman Republic; Stalin was careful not to call himself a Tsar.
In our context, this suggests that there are not four, but sixteen relationships between the %structures of 
key components of succeeding theories %(where the subscript indicates the apparent relationship):
(using subscripts to indicate outward appearances):
%\[
%\begin{array}{cccc}
%=/=&\cong/=&\hookrightarrow/=&\hookleftarrow/=\\
%=/\cong&\cong/\cong&\hookrightarrow/\cong&\hookleftarrow/\cong\\
%=/\hookrightarrow&\cong/\hookrightarrow&\hookrightarrow/\hookrightarrow&\hookleftarrow/\hookrightarrow\\
%=/\hookleftarrow&\cong/\hookleftarrow&\hookrightarrow/\hookleftarrow&\hookleftarrow/\hookleftarrow
%\end{array}
%\]
%\[
%\begin{array}{cccc}
%=_=&\cong_=&\hookrightarrow_=&\hookleftarrow_=\\
%=_{\cong}&\cong_{\cong}&\hookrightarrow_{\cong}&\hookleftarrow_{\cong}\\
%=_{\hookrightarrow}&\cong_{\hookrightarrow}&\hookrightarrow_{\hookrightarrow}&\hookleftarrow_{\hookrightarrow}\\
%=_{\hookleftarrow}&\cong_{\hookleftarrow}&\hookrightarrow_{\hookleftarrow}&\hookleftarrow_{\hookleftarrow}
%\end{array}
%\]
\[
\begin{array}{cccc}
K_{n}=_=K_{n+1}&K_{n}\cong_=K_{n+1}&K_{n}\hookrightarrow_=K_{n+1}&K_{n}\hookleftarrow_=K_{n+1}\\
K_{n}=_{\cong}K_{n+1}&K_{n}\cong_{\cong}K_{n+1}&K_{n}\hookrightarrow_{\cong}K_{n+1}&K_{n}\hookleftarrow_{\cong}K_{n+1}\\
K_{n}=_{\hookrightarrow}K_{n+1}&K_{n}\cong_{\hookrightarrow}K_{n+1}&K_{n}\hookrightarrow_{\hookrightarrow}K_{n+1}&K_{n}\hookleftarrow_{\hookrightarrow}K_{n+1}\\
K_{n}=_{\hookleftarrow}K_{n+1}&K_{n}\cong_{\hookleftarrow}K_{n+1}&K_{n}\hookrightarrow_{\hookleftarrow}K_{n+1}&K_{n}\hookleftarrow_{\hookleftarrow}K_{n+1}
\end{array}
\]
The salutary point here is that we need to be careful in how we specify our terms, lest we misclassify (apparent) revolutions. In particular, the sharing of `lexical taxonomy' had better be more than just lexical, or the innocuous null revolution at the bottom left could count as a case of TI while the stealthy inglorious revolution at the top right does not.

A bold sceptical thesis would be that the rightmost two columns are in practice uninstantiated; that is, all revolutions are glorious, and all appearances to the contrary deceptive.
As we will see in the next section, just such a view has been proposed by some historians of mathematics. %In the remainder of the chapter we will see how well it stands up.

\section{Mathematical Revolutions?}\label{Sec:Revolutions}
The discussion of mathematical revolutions essentially begins with Michael Crowe, who boldly asserts as a `law' that `Revolutions never occur in mathematics' \cite[19]{Crowe75}. Nonetheless, his own subsequent writings are increasingly nuanced: he has moved from denying that there any revolutions in mathematics %\citeyearpar[19]{Crowe75} 
to suggesting that even inglorious revolutions may be possible  (\citeyear[264~f.]{Crowe88}; \citeyear[313]{Crowe92}). % \cite{Crowe88,Crowe92}.
Joseph Dauben has published several articles arguing for the existence of mathematical revolutions \cite{Dauben84,Dauben92,Dauben96}. However, as we shall see, his position is much closer to Crowe's than might be expected.
Indeed, as one commentator, writing more than twenty years ago, has remarked, the literature on mathematical revolutions represents an `authentic theoretical shambles' \cite[193]{Otero96}.
%This does not seem unduly pessimistic. 
There have been two collections of papers on the topic \cite{Gillies92,Ausejo96}.
But the editor of the first notes that each of his authors `has a different theoretical perspective% on the question of revolutions in mathematics
' \cite[8]{Gillies92}. And a contributor to the second pointedly observes that none of these authors make much use of a Kuhnian framework in their other writing on mathematical practice \cite[170]{Corry96}.
In this section I shall attempt to resolve some of this confusion.%
\footnote{For a more extensive account of the debate over revolutions in mathematics, see \cite[107 ff.]{Francois10}.}

%\cite{Corry93,Corry96}.
Many accounts of revolutions in mathematics distinguish two sorts of revolution, usually in terms of the presence or absence of some sort of conceptual continuity. %(or non-preservation).
Hence Crowe distinguishes a `transformational event', in which `an accepted theory is overthrown by another theory, which may be old or new', from a `formational event', in which `an area of science is not transformed, but is \emph{formed}. The discovery that produces this effect is usually new, and by definition overthrows and replaces nothing' \cite[310, citing his own earlier work]{Crowe92}.
%Subsequent commentators (for instance, \cite[123~f.]{Crowe67}; \cite[5]{Gillies92}) have argued that, although revolutions of this character do occur, there can also be revolutionary change in which all the concepts are retained, albeit with a transformed character. This distinction is a familiar one in political history, where the revolutionary metaphor originates. We may distinguish between the Russian Revolution of 1917, in which the whole constitution was abandoned and replaced by something radically different, with different constituent parts, and the Glorious Revolution of 1688 in Britain, in which all the principal constitutional constituents, the crown, parliament, and so forth, were retained, although their character and relative significance changed dramatically.
Dauben likens mathematical revolutions to the Glorious Revolution of 1688, in the persistence of the `old order', albeit `under different terms, in radically altered or expanded contexts' \cite[52]{Dauben84}.
This political analogy is echoed by Donald Gillies, who frames the distinction as between `Franco-British' and `Russian' revolutions: %in similar terms to our contrast of `glorious' and `inglorious' revolutions \cite[5]{Gillies92}. %Our terminology may exhibit unabashed persuasive definition, but it sidesteps the historical exegesis prompted by Gillies: why is
in the former a `previously existing entity persists' %, but experiences 
through `a considerable loss of importance'; 
in the latter the `previously existing entity' is `overthrown and irrevocably discarded' 
\cite[5]{Gillies92}. 
Gillies's choice of terminology provokes the distracting historiographical question of why the French revolution should be more like the British than the Russian. %? What are the start and end points of each revolution?
After all, France and Russia both ended up as republics, whereas Britain did not. The answer seems to be that Gillies stops the clock at some point in the reign of Louis-Philippe.
%Which is why Gillies \citeyear[5]{Gillies92} thinks the French Revolution was glorious, since he includes the 1815 restoration of Louis XVIII within its scope. However, this indicates the instability of assessments made on such a large scale, since there seems no good reason why he should not have gone further still and included the overthrow of Louis-Philippe in 1848, which makes the whole affair inglorious.
My own use of `glorious' revolutions, inspired by Dauben's usage but not intended to be given any specific historical reading, at least has the merit of sidestepping such musings.
%In Crowe's terminology glorious and (tacitly) paraglorious revolutions are `formational events' and inglorious revolutions are `transformational events.'
More importantly, we may notice that these distinctions do not necessarily coincide---Crowe's formational events seem closer to paraglorious than glorious revolutions, for example---and that, although all of them are binary, none of them seems to be exhaustive. So a more fine-grained distinction may be a source of clarity.

However the distinction %between glorious and inglorious revolutions 
is drawn, most of its framers agree that only glorious revolutions are possible in mathematics:
`One important consequence, in fact, of the insistence on self-consistency within mathematics is that its advance is necessarily cumulative. New theories cannot displace the old' \cite[62]{Dauben84};
`In science both Russian and Franco-British revolutions occur. In mathematics, revolutions do occur but they are always of Franco-British type' \cite[6]{Gillies92};
`revolutions do occur in mathematics, but are confined entirely to the metamathematical component of the community's shared background' \cite[223]{Dunmore92}. 
In this regard they concur with the opinion of many mathematicians that their discipline is cumulative. Crowe finds such sentiments expressed by Fourier in 1822, Hankel in 1869, and Truesdell in 1968 \cite[19]{Crowe75}. 
Indeed, celebrated mathematicians are still saying as much: `central contributions have been lasting, one does not supersede another, it enlarges it' \cite[25]{Langlands13}.
As Crowe poetically summarizes the conventional view,
`Scattered over the landscape of the past of mathematics are numerous citadels, once proudly erected, but which, although never attacked, are now left unoccupied by active mathematicians' \cite[263]{Crowe88}.
Nonetheless, there are exceptions to this trend. Michael Harris notes Kronecker in 1891 observing that `in this respect mathematics is no different from the natural sciences: new phenomena overturn the old hypotheses and put others in their place' and Siegel in 1964 characterizing work revisionary of his own as `a {pig} broken into a beautiful {garden} and rooting up all flowers and trees' \cite[4]{Harris15}.
%The ground for denying that inglorious revolutions occur in mathematics is that the discipline is cumulative in a way that empirical science is not: both disciplines discard old material, but mathematicians never really throw it away. 
%Quaternions or conic sections may be of no greater interest to the modern mathematician than phlogiston or caloric are to the modern physicist, but their legitimacy is not disputed.

Bruce Pourciau complains that the `Crowe--Dauben debate' is actually a `Crowe--Dauben consensus', viz.: {Kuhnian revolutions are inherently impossible in mathematics}' \cite[301]{Pourciau00}.
For Pourciau a revolution is Kuhnian or
`\emph{noncumulative} whenever some true statements of the old conception have no translations (faithful to the original meaning) which are true statements in the new conception' \cite[301]{Pourciau00}.
Pourciau argues that Brouwerian intuitionism %, which, if adopted, would have required wholesale revision of results treated as certain by prior mathematicians, 
is a (failed) Kuhnian revolution.
Certainly, if adopted, this would have required wholesale revision of results treated as certain by prior mathematicians, thereby meeting the strictest definition of Kuhnian revolution. Its usefulness as an example might be somewhat compromised by the fact that it never actually happened, but it was seriously proposed and still commands some support.
However, Pourciau may be overestimating the difficulty in supplying examples of Kuhnian revolutions in mathematics in two ways: one of scale, one of chronology.
Firstly, although the standard example of a Kuhnian revolution in natural science is the Copernican revolution, an epochal upending of an all-encompassing worldview, it is a mistake to suppose that all Kuhnian revolutions need be so drastic.
Stephen Toulmin once complained that Kuhn had surreptitiously revised his position to admit `small-scale ``micro-revolutions''' \cite[47]{Toulmin70}.
Kuhn strenuously rejected this imputation: `My concern \dots\ has been throughout what Toulmin now takes it to have become: a little studied type of conceptual change which occurs frequently in science and is fundamental to its advance' \cite[249 f.]{Kuhn70}.
He subsequently characterized a paradigm as `what the members of a scientific community and they alone share' %\cite[460]{Kuhn77}, 
where such communities may comprise `perhaps 100 members, sometimes significantly fewer' \cite[460; 462]{Kuhn77}.
Happily enough, this coincides with an influential estimate of the size of mathematical research communities: 
%`The objects which our mathematician studies were unknown before the twentieth century; most likely, they were unknown even thirty years ago. Today they are the chief interest in life for 
`a few dozen (at most a few hundred)% of his comrades
' \cite[35]{Davis80}.
%Conversely, Dauben may be underestimating the revolutionary character of some of his examples.
%\cite{Quinn12,Buldt16}. \cite{Segura16}.
Secondly, as observed above, paraglorious and inglorious revolutions are essentially symmetrical; they differ only in the chronological sequence of the contrasting theories. 
Paraglorious revolutions are cumulative, but they exhibit a conceptual discontinuity formally identical to that exhibited by inglorious revolutions.
Chronological sequence alone does not seem to be a principled basis on which to discount paraglorious revolutions as Kuhnian.% (whether or not they would count as revolutions for Kuhn).%
\footnote{I make no claim as to Kuhn's own view on this issue, although I note that he does refer to historians experiencing revolutions by `moving through time in a direction opposite to scientists' \cite[57]{Kuhn83}, which at least suggests an openness to temporal symmetry.}

%`From an \emph{epistemological} point of view, not only do the various paradigms \dots\ in mathematics \dots\ not cancel each other out, in fact, they benefit from \emph{mixed fusions}' \cite[301]{Zalamea12}.

%`relate the quasicrystal revolution to the non-Euclidean geometry revolution' \cite{Ashkenazi14}.

%\cite{Kanamori12}.
%\cite{Weber13}.
%\cite{Zeilberger01}.

%`a Fregean revolution in logic' \cite{Anellis11}.
%`a Toulminian revolution in logic' \cite{Castaneda60}.
%`a current revolution in logic' (linear logic) \cite{Fraser10}. \cite{Castellana11}.

Three broad strategies for the identification of Kuhnian (or nonglorious) revolutions in mathematics arise from this discussion. Firstly, we may look directly for inglorious revolutions: conceptual shifts within mathematics in which key components have been lost. 
Secondly, we may look for paraglorious revolutions: conceptual shifts within mathematics in which key components have been gained. 
Thirdly, we may look for sorites-like sequences of glorious (or paraglorious) revolutions which exhibit non-transitivity of glory, that is which are collectively inglorious.
The search is complicated by several factors. In particular, it is not easy to determine whether a given episode is revolutionary; nor is it easy to determine what type of revolution a given revolutionary episode exemplifies. Hence some of the same examples might be claimed as successes for more than one of these search strategies. An analysis of even a single case study thorough enough to settle all of these issues would be beyond the scope of this chapter. However, in the following sections I will discuss several putative mathematical revolutions in what I hope to be at least enough detail to indicate %some of the directions a more thorough approach might pursue.
the prospects for these strategies.

\section{$\mathbb{Q}\rightarrow\mathbb{R}\rightarrow\mathbb{C}$}%Original Incommensurability}
The most obvious example of incommensurability in mathematics must be incommensurability itself! The concept is, of course, originally a mathematical one, credited to the very earliest Greek geometers, the Pythagoreans. Thomas Heath, in his edition of Euclid, quotes a scholium on the first proposition of Book X, attributed to Proclus:
`They called all magnitudes measurable by the same measure commensurable, but those which are not subject to the same measure incommensurable' \cite[684]{Heath06}.
Specifically, the Pythagoreans discovered that $\sqrt{2}$ was incommensurable with the natural numbers, that is, it cannot be expressed as a ratio of natural numbers or, as we would say, as a rational number. 
The discovery was credited to one Hippasus of Metapontum, who is reputed to have drowned in a shipwreck.
As the historian of mathematics Kurt von Fritz observes, 
\begin{quotation}
The discovery of incommensurability must have made an enormous impression in Pythagorean circles because it destroyed with one stroke the belief that everything could be expressed in integers, on which the whole Pythagorean philosophy up to then had been based. This impression is clearly reflected in those legends which say that Hippasus was punished by the gods for having made public his terrible discovery \cite[260]{Fritz45}.
\end{quotation}
Even in ancient times, the allegorical aspects of this story were already apparent,
%`the first of the Pythagoreans who made public the investigation of these matters perished in a shipwreck' 
`hinting that everything irrational and formless is properly concealed, and, if any soul should rashly invade this region of life and lay it open, it would be carried away into the sea of becoming and be overwhelmed by its unresting currents', as Proclus puts it \cite[684]{Heath06}.

For our purposes, the crucial point in this narrative is that the change initiated by Hippasus was revisionary of earlier mathematics: it %. It was not a conservative extension of the existing concept of number. Rather, it
%`The extension of the theory of proportion to incommensurables 
`required an entirely new concept of ratio and proportion and a new criterion to determine whether two pairs of magnitudes which are incommensurable with one another have the same [ratio]' %\emph{logos}' 
\cite[262]{Fritz45}.
The completion of this task by later mathematicians, notably Theaetetus and Eudoxus, is one of the great achievements of Greek mathematics, and plausibly a major driver of its early development of the concept of rigorous proof.
For Dauben it is one of the best examples of a mathematical revolution. He stresses that the
%`The old concept of number, although the word was retained, was gone, and in its place, numbers included irrationals as well. This 
`transformation of the concept of number \dots\ %, however, 
entailed more than just extending the old concept of number by adding on the irrationals---the entire concept of number was inherently changed, transmuted as it were, from a world-view in which integers alone were numbers, to a view of number that was eventually related to the completeness of the entire system of real numbers' \cite[57]{Dauben84}.

In my terminology, this is clearly not a glorious revolution, since a literally incommensurable concept has been added. So it is at least a paraglorious revolution. Might we go further and identify it as also inglorious? On the one hand something has certainly been lost: the `world-view in which integers alone were numbers', for a start. On the other hand, world-views are not part of the subject matter of mathematics. Hence Caroline Dunmore identifies this shift as `the first great meta-level revolution in the development of mathematics' \cite[215]{Dunmore92}. On Dunmore's account, object level revolutions in mathematics are always glorious, but they are always accompanied by inglorious revolutions in the meta-level, that is in the philosophical or methodological presuppositions \cite[225]{Dunmore92}. 
The rational numbers are still an object of mathematical enquiry %---they haven't gone away---
and the Pythagorean results about their comparison still hold. Nonetheless, it is highly misleading to conceive of the real numbers as a conservative extension of the rationals. The real numbers are constructed on a quite different basis, but in such a way that a subset isomorphic to the rationals may be identified. 

Strictly speaking, the taxonomic incommensurability between mathematics defined over $\mathbb{Q}$ and mathematics defined over $\mathbb{R}$ runs both ways. Clearly, real mathematics cannot be done with rationals alone, but techniques that work over $\mathbb{Q}$ fail over $\mathbb{R}$.
So a revolution from the mathematics of $\mathbb{Q}$ to the mathematics of $\mathbb{R}$ would be paraglorious and inglorious. 
However, the actual revolution %was messier, and 
could also be described as a shift from the mathematics of $\mathbb{Q}$ to (eventually) the mathematics of $\mathbb{Q}$ and $\mathbb{R}$, understood as separate projects. That shift would be strictly paraglorious---assuming that the mathematics of $\mathbb{Q}$ has been preserved, and not just reconstructed.
The same issue arises with %other types of number, 
supersets of the reals, whether well-established, such as the complex numbers,
or more contentious, such as the hyperreals, which include infinitesimals \cite{Bair13}. %\cite{Blaszczyk13}. \cite{Katz13}, \cite{Katz13a}.
The underlying issue of cross-sortal identity is a known problem for a wide range range of philosophies of mathematics \cite[124]{Cook05a}.
Textbook presentations of the foundations of mathematics are obliged to address cross-sortal identity, which they do in a variety of ways, often at odds with mathematical practice \cite[for a careful discussion, see][180 ff.]{Ganesalingam13a}.
It is also important to note that retrofitting a new foundation onto existing mathematics is not confined to number systems. Indeed it has been a major feature of mathematical research since the nineteenth century---and it is precisely what Siegel was complaining about as `rooting up all flowers and trees' %in the remark quoted %in \S\ref{Sec:Revolutions} above 
\cite{Lang94}.
It is a deep question whether such moves can be understood as merely adding `a new storey to the old structure' \cite[19, quoting Hankel]{Crowe75}. To pursue the architectural metaphor, they might be better characterized as `fa\c{c}ading', whereby the front elevation of an otherwise demolished building is incorporated into its successor. I shan't settle that question here, but we may observe that these shifts are at least paraglorious and perhaps inglorious.

\section{Irony for Mathematicians}\label{Sec:Irony}
One way of approaching the issue of inglorious revolution in mathematics is through a related question: when do mathematicians say things that are not so?
One prospect might be the assumptions of reductio proofs.
In a recent essay, the mathematician Timothy Gowers briefly considers, but ultimately rejects, 
%`a suggestion %was made \dots\ 
the intriguing idea `that proofs by contradiction are the mathematician's version of irony' \cite[224]{Gowers12}. %I'm not sure that I agree with that: 
He objects that `when we give a proof by contradiction, we make it very clear that we are discussing a counterfactual, so our words \emph{are} intended to be taken at face value' \citetext{op.\ cit.}. %[224]{Gowers12}.
Perhaps more tellingly, we might frame this objection as saying that a proposition assumed for the purposes of proof by contradiction is presented as the antecedent of a conditional: `If $P$ were the case, then \dots\ a contradiction would follow. So, not $P$.'
Nonetheless, at least for the duration of the ellipsis, the mathematician proceeds as though $P$ were being seriously entertained. An inattentive reader who began reading a proof part way through would not necessarily be able to tell which proposition the mathematician intended to show to be false---or even that any of them were presented with this intent.

Another possibility might be unproven conjectures upon which mathematicians sometimes rely, when exploring their consequences. The mathematician Barry Mazur talks of `architectural conjectures' that `play the role of ``joists'' and ``supporting beams'' for some larger mathematical structure yet to be made' \cite[199]{Mazur97}.
The formulation of such conjectures
%`These conjectures are \emph{expected} to turn out to be true, as, of course, are all conjectures; their formulation 
\begin{quotation}
is often a way of ``formally'' packaging, or at least acknowledging, an otherwise shapeless body of mathematical experience that points to their truth. \dots\ %From these conjectures, implications may be perfectly rigorously made. Best, if the conjectures are, loosely speaking, ``testable'', or ``falsifiable'' in the sense that they imply a stream of particular, numerical perhaps, predictions many of which may be directly checked. But these conjectures are \emph{architectural} in that they play the role of ``joists'' and ``supporting beams'' for some larger mathematical structure yet to be made. 
These conjectures sometimes round out a field by being clear, general (but not yet proved) statements enabling one to understand where a certain amount of on-going, perhaps fragmentary, specialized work is headed; they provide a focus \citetext{op.\ cit.}. %\cite[199]{Mazur97}.
\end{quotation}
Architectural conjectures seldom arise alone; they often comprise elaborate networks of interlinked conjectures that present the outline of what is hoped to be many years of fruitful work. One of the best known such networks of conjectures in contemporary mathematics is
the {Langlands programme}, `an extensive web of conjectures by which number theory, algebra, and analysis are interrelated in a precise manner, eliminating the official divisions between the subdisciplines' \cite[180]{Zalamea12}. This has been enormously influential, guiding the work of scores of mathematicians who have confirmed some---but by no means all---of its key conjectures.

As with reductio hypotheses, conjectures are strictly to be understood as the antecedents of conditionals. Mathematicians should not be seen as \emph{asserting} them until they have actually been proven. Nonetheless, as with reductio hypotheses, they are presented in apparent earnest, and their implications investigated with all due rigour: they `are \emph{expected} to turn out to be true, as, of course, are all conjectures' \cite[199]{Mazur97}.
So, naively, it may seem as though neither sort of hypothesis is much use for present purposes, since mathematicians seem to have an uncanny knack of only assuming for proof by contradiction things that are false and only assuming as conjectures things that are true (even if not yet proven).
This would be a profound misperception:
attempted reductios sometimes founder on the truth of the hypothesis 
(famously so in the accidental discovery of non-Euclidean geometry)
and sometimes substantial effort is devoted to exploring the consequences of false conjectures.
In the next section we will encounter an example of the latter.

%Three grades of mathematical irony
%\begin{enumerate}
%\item Proofs by contradiction
%\item Quandaries
%\item False conjectures
%\end{enumerate}
%In each case a mathematician asserts that $p$ where $p$ is false. 
%In (1) the mathematician is justified in believing $\neg p$;
%in (2) the mathematician is justified in believing neither $p$ nor $\neg p$;
%in (3) the mathematician is justified in believing $p$.

%`How should one regard a proof if it relies on the Riemann hypothesis? One could simply say that the proof establishes that such and such a result is implied by the Riemann hypothesis and leave it at that. But most mathematicians take a different attitude. They believe the Riemann hypothesis, and believe that it will one day be proved. So, they believe all its consequences as well, even if they feel more secure about results that can be proved unconditionally' \cite[68]{Gowers08}.
%The `{Langlands program} is a collection of conjectures, due to Robert Langlands, that relate number theory to representation theory \dots\ Between them, these conjectures generalize, unify, and explain large numbers of other conjectures and results. For example, the Shimura-Taniyama-Weil conjecture, which was central to Andrew Wiles's proof of Fermat's Last Theorem, forms one small part of the {Langlands program}' \cite[69]{Gowers08}.

\section{The World Without End Hypothesis}
In 2016 Michael Hill, Michael Hopkins, and Douglas Ravenel published an article in %\emph{Annals of Mathematics}, 
one of the most prestigious journals in mathematics, with the following abstract:
`We show that the Kervaire invariant one elements $\theta_j \in \pi_{2^{j+1}-2}S^{0}$ exist only for $j\leq6$. By Browder's Theorem, this means that smooth framed manifolds of Kervaire invariant one exist only in dimensions 2, 6, 14, 30, 62, and possibly 126. Except for dimension 126 this resolves a longstanding problem in algebraic topology' \cite[1]{Hill16}.
This was the final, fully vetted version of a result that they had announced seven years earlier. As they summarize their result in a preliminary expository article, they showed that
`certain long sought hypothetical maps [the $\theta_j$ for $j\geq7$] between high dimensional spheres do not exist' \cite[32]{Hill09}.
This outcome was a surprising one, not just because of the technical depth of the work required but because many experts in the area had long expected the opposite result:
`The problem solved by our theorem is nearly 50 years old. There were several unsuccessful attempts to solve it in the 1970s. They were all aimed at proving the opposite of what we have proved' \cite[32]{Hill09}.
The hypothesis that the sought after maps all exist came to be known as the World Without End Hypothesis; the contradictory hypothesis that the $\theta_j$ only exist for small $j$ was known as the Doomsday Hypothesis. Hill, Hopkins, and Ravenel proved the Doomsday Hypothesis, and thereby disproved the World Without End Hypothesis.

%Indeed, a body of work had been constructed on that
While not remotely on the scale of the Langlands Programme, the World Without End Hypothesis was not just a single assertion, but the basis for a whole system of `architectural conjectures': the new proof demolished
%`Researchers developed 
`what Ravenel calls an entire ``cosmology" of conjectures% based on the assumption that manifolds with Arf-Kervaire invariant equal to 1 exist in all dimensions of the form $2^{n}-2$
' \cite[374]{Klarreich11}.
The triumph of the Doomsday Hypothesis undercut a growing sense of understanding provided by the World Without End Hypothesis. As the reviewer of Hill, Hopkins, and Ravenel's paper in \emph{Mathematical Reviews} comments, the World Without End Hypothesis
`was so compelling that many believed the $\theta_j$ must exist; now that we know they don't, the behavior of the EHP sequence is much more mysterious. In particular, Mahowald's $\eta_j$-elements \dots\ %[M. E. Mahowald, Topology 16 (1977), no. 3, 249Ð256; MR0445498] 
now appear entirely anomalous' \cite{Goerss16}.
The author of a book surveying some of the techniques developed in pursuit of the World Without End Hypothesis published %as the balance of likelihood seemed to be tipping towards the Doomsday Hypothesis 
shortly before Hill, Hopkins, and Ravenel's announcement concluded his preface as follows:
`In the light of the above conjecture [the Doomsday Hypothesis] and the failure over fifty years to construct framed manifolds of Arf-Kervaire invariant one this might turn out to be a book about things which do not exist' \cite[ix]{Snaith09}.

%Nonetheless, as far as that community is concerned, this Small as it may be, t
The %confirmation of the Doomsday Hypothesis and the 
collapse of the World Without End Hypothesis seems to be an inglorious revolution. 
It is clearly a case of referential discontinuity, since a whole class of objects which were key to the old theory have been shown not to exist. 
It might be objected that there is no taxonomic incommensurability, since the conjectures are still perfectly intelligible, despite now being known to be false. 
Indeed we have seen that Mizrahi %rightly 
argues that referential discontinuity does not entail conceptual incompatibility, and is thereby %does not support taxonomic incommensurability (see \S\ref{Sec:Mizrahi} above). In other words, referential discontinuity is 
insufficient for %conceptual incompatibility 
taxonomic incommensurability (see \S\ref{Sec:Mizrahi} above). So what else is required?
A natural candidate would be Kuhn loss: reduction in actual or potential explanatory power. %loss of `some actual and much potential explanatory power' \cite[107]{Kuhn62}. 
Kuhn loss does not seem to feature in the examples Mizrahi adduces of referential discontinuity without conceptual incompatibility \cite[365 ff.]{Mizrahi14b}, but it is exhibited here: the understanding that had seemed to follow from the World Without End Hypothesis has been lost. 
Specifically, comments such as Paul Goerss's reveal the `mysterious' and `anomalous' nature of some of the surviving results, now that the architectural conjectures which mathematicians had relied on to understand them have been falsified. This is what we should expect from a reduction in explanatory power.
%Not only was a widely believed conjecture shown to be false, and other conjectures conditional upon it undercut, but the understanding that had seemed to follow from the World Without End Hypothesis has been lost. 
%But that would overlook, firstly, comments such as Paul Goerss's, quoted above, about the `mysterious' and `anomalous' nature of some of the surviving results, now that the architectural conjectures which mathematicians had relied on to understand them have been falsified, and, secondly, that the new taxonomy certainly omits  some things: the %$\theta_j$ for $j\geq7$.
%subjects of the `book about things which do not exist'.

If the collapse of the World Without End Hypothesis is a revolution, it is certainly a small-scale one. %This topic was 
The constituency of homotopy theorists for whom this was an active research area would be at the low end of Kuhn's `perhaps 100 members, sometimes significantly fewer' comprising a scientific community \cite[462]{Kuhn77}.
However, this makes this revolution much more likely to be typical of mathematical revolutions. Large-scale revolutions are necessarily extremely rare; so much so that their instances may well be \emph{sui generis}. Conversely, small-scale revolutions are much better placed to support generalizations: there are plenty of other failed architectural conjectures that could be explored in similar detail \cite[see, for example, some of the `cautionary tales' in][7 f.]{Jaffe93}.

\section{Inter-Universal Teichm\"uller Theory}
One of the more widely discussed open problems in modern number theory is the $abc$ conjecture. Like many such problems it can be stated quite simply, but defies simple solution. Here is what it says:
\theoremstyle{definition} \newtheorem*{dfn}{The $abc$ conjecture}
\begin{dfn}
%$abc$ conjecture. %(Oesterl\'e--Masser, 1985)} 
For every $\varepsilon > 0$, there are only finitely many triples $(a, b, c)$ of coprime positive integers where $a+b = c$, such that $c > d^{1+\varepsilon}$, where $d$ denotes the radical of $abc$ (the product of its distinct prime factors). %of $abc$.
\end{dfn}
For example, try $a=15$ and $b=28$. These are coprime, but $c=43$ and $d=2\times3\times5\times7\times43=9030\gg43$. So $(15,28,43)$ is not one of the specified triples (for any $\varepsilon$). On the other hand, let $a=1$ and $b=63$. Then we have $c=64$ and $d=2\times3\times7=42<64$. So $(1,63,64)$ is such a triple (at least for values of $\varepsilon < .11269$). 

In a series of preprints appearing on his website in 2012, the respected mathematician Shinichi Mochizuki claimed to have a proof of the $abc$ conjecture. However, Mochizuki's claimed proof %of the $abc$ conjecture 
introduced so many new techniques and concepts that other leading mathematicians in the field described it as like `reading a paper from the future, or from outer space' and as `very, very weird' \cite[cited in][]{Chen13}. 
%If Mochizuki has proved the result, it is propositionally (but presumably not doxastically\footnote{Mochizuki would not be doxastically justified if, as seems plausible, ``acting in an epistemically responsible manner'' includes successfully explaining the proof to others and he has not yet done this.}) 
%justified for him, but it's not yet propositionally justified for other mathematicians. %(Of course, if he hasn't proved it, it isn't justified for anybody either way.)
The scale of the proof (more than 500 pages) and its sheer incomprehensibility, even by the standards of cutting-edge research mathematics, have so far stalled all attempts at the normal processes of confirmation and acceptance that transform a proof claim into an established proof.
Although a handful of other mathematicians now profess to understand Mochizuki's work, they have had little success sharing that understanding more widely. One anonymous mathematician, quoted in \emph{Nature}, summed up the problem:
\begin{quotation}
%`But so far, the few who have understood the work have struggled to explain it to anyone else. 
``Everybody who I'm aware of who's come close to this stuff is quite reasonable, but afterwards they become incapable of communicating it''\dots\ %,''says one mathematician who did not want his name to be mentioned. 
The situation, he says, reminds him of the \emph{Monty Python} skit about a writer who jots down the world's funniest joke. Anyone who reads it dies from laughing and can never relate it to anyone else \cite[181]{Castelvecchi15}.
\end{quotation}

%`Mochizuki's four IUT theory papers \dots\ %, which 
%were the subject of the last two days of the conference. The job of explaining those papers fell to Chung Pang Mok of Purdue University and Yuichiro Hoshi and Go Yamashita, both colleagues of Mochizuki's at the Research Institute for Mathematical Sciences at Kyoto University. The three are among a small handful of people who have devoted intense effort to understanding Mochizuki's IUT theory. By all accounts, their talks were impossible to follow' \cite{Hartnett15}.

Mochizuki calls his work inter-universal Teichm\"uller theory, or IUTeich.
He has reflected on the verification process of IUTeich in a pair of papers that comprise a valuable resource for philosophers of mathematical practice \cite{Mochizuki13,Mochizuki14}.%
\footnote{For a very different application of these papers to the philosophy of mathematical practice, see \cite[187 ff.]{Tanswell16a}.}
Mochizuki warns that `any attempt to study IUTeich under the expectation that the \emph{essential thrust} of IUTeich will proceed via a similar pattern of argument to existing mathematical theories is likely to end in failure' \cite[5; all emphases Mochizuki's]{Mochizuki13}. Even Teichm\"uller theory itself is only an indirect inspiration.
Nonetheless, IUTeich does echo Teichm\"uller theory in at least one respect---% 
%`Teichm\"uller's (1939, 1943) papers 
the papers in which Oswald Teichm\"uller laid out his theory were not immediately accepted by the mathematical community either:
%`were first considered by the other mathematicians as a program and not as finished papers. \dots\ %According to Abikoff, the fact that several ideas were sketched is only consistent with the tradition of the journal in which the article was published. Abikoff (1986) writes: ÒThe tradition of Deutsche Mathematik is one of heuristic argument and contempt for formal proof. Busemann notes that TeichmŸller manifested those traits early in his career but when pressed could offer a formal proof.Ó 
`It was after several years of hard work by several mathematicians that all the arguments in these papers were considered as being sound' \cite[128]{Ji13}.
So a lengthy gap %between initial announcement and 
before final community acceptance is not unusual in itself. (Note also the seven years between initial announcement and final publication of Hill, Hopkins, and Ravenel's work.)
What is unusual in Mochizuki's case is that the mathematical community appears to be completely stumped.

The trouble arises from both the scale and the nature of the task required of mathematicians who wish to come to terms with Mochizuki's work. He suggests, perhaps optimistically, that
`it is quite \emph{possible} to achieve a \emph{reasonably rigorous understanding} of the theory within a period of \emph{a little less than half a year}' \cite[4]{Mochizuki13}.
But this is still a substantial investment of time. Mochizuki also notes that his work is essentially independent of the Langlands Programme (discussed in \S\ref{Sec:Irony}) \cite[10]{Mochizuki14}. Since this has guided so much recent work in number theory, many of the individuals most interested in the $abc$ conjecture have a background that does not particularly suit them for tackling IUTeich, and should not necessarily expect to acquire techniques that would further their own projects from the six months or more of concentrated intellectual effort required.
Indeed, Mochizuki stresses the incompatibility of the ideas behind IUTeich and the ideas most number theorists are familiar with:
`the most essential stumbling block lies not so much in the need for the \emph{acquisition of new knowledge}, but rather in the need for researchers \dots\ %(i.e., who encounter substantial difficulties in their study of IUTeich and related topics)
to \emph{deactivate the thought patterns} that they have installed in their brains and taken for granted for so many years and then to start afresh%, that is to say, to \emph{revert to a mindset} that relies only on \emph{primitive logical reasoning}, in the style of a student or a novice to a subject
' \cite[11 f.]{Mochizuki14}.
He complains that
\begin{quotation}
%Typically, 
when a researcher with a solid track record in mathematical research decides to read a mathematical paper, \dots\ %unlike the case with students or novices who take the time to \emph{study step by step from the rudiments of a subject}, 
such a researcher will attempt to digest the content of the paper in as efficient a way as is possible, by \emph{scanning} the paper for important terms and theorems so that the researcher may apply his/her vast store of expertise and deep understanding of the subject to determine just \emph{which} of those topics of the subject that, from point of view of the researcher, have already been ``digested'' and ``well understood'' \emph{play a key role in the paper}. \dots\ 
%Put another way, this amounts to the sort of ``\emph{occasional nibbling}'' that Yamashita warned of during his lecture series at Kyushu University (cf. (2) above). 
Of course, in the case of IUTeich, a researcher who already possesses a deep understanding, as well as a solid track record in mathematical research, concerning such topics as absolute anabelian geometry, the rigidity properties of the \'etale theta function, and Hodge-Arakelov theory, may indeed find such ``occasional nibbling'' to be more than sufficient to attain a quite genuine understanding of IUTeich. In fact, however, for better or worse, \emph{no such researcher exists} (other than myself) at the present time
\cite[8 f.]{Mochizuki14}.
\end{quotation}
This is an eloquent description of conceptual incompatibility. Mochizuki is not just saying that no one is capable of understanding IUTeich; on the contrary, he is confident that the material is well within the grasp of competent research mathematicians. But he stresses that, if they are to understand IUTeich, they must sat aside their existing conceptual frameworks and build a new one from scratch.
The contrast which Mochizuki draws echoes that Kuhn draws between translation as opposed to interpretation and language acquisition; incommensurability being a bar to the former, but not the latter \cite[53]{Kuhn83}.

If the revolution which Mochizuki has so far failed to ignite succeeds, it would appear to be strictly paraglorious in nature. He does not wish to overturn anything; rather he wishes to comprehensively supplement the existing apparatus of number theory. If IUTeich is correct, it will represent a substantial leap forward in mathematics. Mochizuki's problem is that he's trying to do it all in one go. Conversely, if an irreparable flaw is found in Mochizuki's reasoning, and IUTeich collapses, then its fall would be an inglorious revolution.
In either case, IUTeich would exhibit similar conceptual incompatibility with mainstream %representation-theoretic 
number theory.% as pursued within the Langlands Programme.

%`Certain researchers believe that every essential phenomenon in number theory may in fact be reduced to some aspect of the \emph{representation-theoretic} approach exemplified by the \emph{Langlands program}. On the other hand, the fundamental ideas of IUTeich are not based on this sort of representation-theoretic approach' \cite[10]{Mochizuki14}.
%The {Langlands program} is `an extensive web of conjectures by which number theory, algebra, and analysis are interrelated in a precise manner, eliminating the official divisions between the subdisciplines' \cite[180]{Zalamea12}.
%`By contrast [with Wiles's fruitful proof of Fermat's Last Theorem], it is by no means clear that such extensions and generalizations will be possible in the case of IUTeich' \cite[10]{Mochizuki14}.

\section{Classical $\rightarrow$ Modern $\rightarrow$ Contemporary}
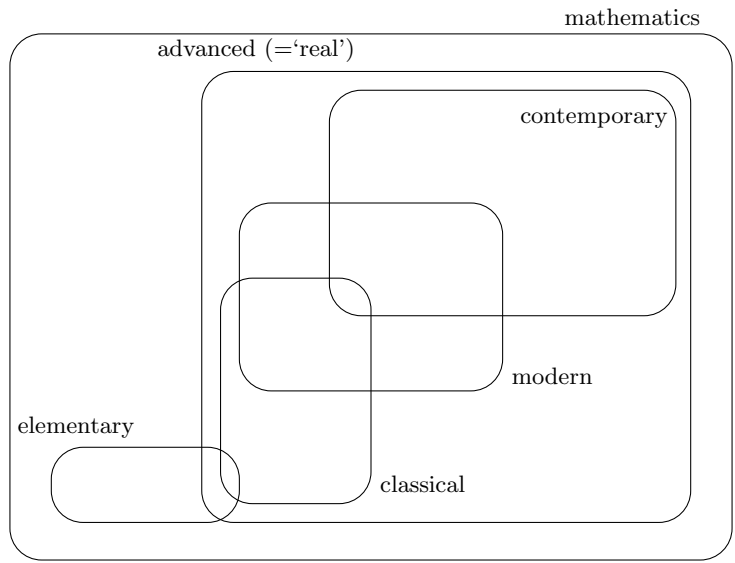
\begin{figure}[htbp]
\small
\begin{center}
\begin{tikzpicture}[every node/.style={draw,rounded corners=12pt}] 
\node at (0,0) [minimum height=7cm,minimum width=9.6cm,label=55:mathematics] {};
\node at (-3,-2.5) [minimum height=1cm,minimum width=2.5cm,label=93:elementary] {};
\node at (1,0) [minimum height=6cm,minimum width=6.5cm,label=110:{advanced (=`real')}] {};
\node at (-1,-1.25) [minimum height=3cm,minimum width=2cm,label=below right:classical] {};
\node at (0,0) [minimum height=2.5cm,minimum width=3.5cm,label=-25:modern] {};
\node at (1.75,1.25) [minimum height=3cm,minimum width=4.5cm,align=right%,label distance=-1cm,label=center:contemporary
] {\hspace{66pt} contemporary\\\\~\\\\\\\\};
\end{tikzpicture}
\caption{Correlations between the areas of mathematics: elementary, advanced, classical, modern, contemporary \protect\cite[after][26; used with permission]{Zalamea12}}
\label{Zalamea}
\end{center}
\end{figure}

As a final example, I wish to shift focus from some ultimately quite small-scale revolutions (albeit ones that have attracted a fair bit of publicity) to the discipline of mathematics as a whole.
The mathematician Fernando Zalamea offers the following useful periodization of research in his discipline:
\begin{description}
\item[Classical mathematics (midseventeenth to midnineteenth centuries):] sophisticated use of the infinite (Pascal, Leibniz, Euler, Gauss);
\item[Modern mathematics (midnineteenth to midtwentieth centuries):] sophisticated use of structural and qualitative properties (Galois, Riemann, Hilbert);
\item[Contemporary mathematics (midtwentieth century to present):] sophisticated use of the properties of transference, reflection and gluing (Grothendieck, Serre, Shelah).
 \cite[27]{Zalamea12}
\end{description}
%Furthermore, t
The scale of mathematical research grows with each generation, such that each period attacks a broader front and produces more results than its predecessors. Zalamea graphically represents this process in the diagram reproduced as Fig.~\ref{Zalamea}.
The moral for the philosopher of mathematical practice is a striking one: classical and modern mathematics may be familiar enough from school or undergraduate study, but contemporary mathematics almost certainly is not. Even philosophers %with a tertiary level education in mathematics 
who take pains to reflect more than just elementary and foundational work may well be quite out of touch with the conceptual underpinnings of mathematical research conducted in their own lifetimes.
Hence, as Zalamea complains, the large scale conceptual shift from modern to contemporary mathematics has gone largely unremarked by philosophers. I would contend that this shift may be understood as revolutionary. 
Certainly the `properties of transference, reflection and gluing' would be impossible to articulate with only the conceptual resources available to Galois, Riemann, or Hilbert (let alone Pascal, Leibniz, Euler, or Gauss). Insofar as these are key components of contemporary mathematics, their acquisition is at least paraglorious.
Furthermore, contemporary mathematics has undergone a sorites-like sequence of paraglorious revolutions of such daunting scope that the key components preserved throughout the sequence have a drastically diminished role in the new era. So much so indeed, that the whole transition might best be characterized as inglorious. 

\section{Conclusion}\label{Sec:Con}
To take stock, we have seen four case studies exemplifying different classes of putative revolution in mathematics: the shift from rational to real numbers (and other cases of foundational retrofitting); shifts occasioned by the collapse of an architectural conjecture, such as the World Without End Hypothesis; shifts resulting from a rapid advance, such as IUTeich; and the collective large-scale shift that has transformed recent mathematics.
The first of these is at least paraglorious and perhaps also inglorious.
The second seems to be strictly inglorious whereas the third is strictly paraglorious (if successful; if unsuccessful, it would be another failed architectural conjecture). Both of these examples are comparatively small-scale %, which suggests that similar shifts may be found
and might be seen as exemplary of similar shifts in other areas of mathematics.
Lastly, the shift from modern to contemporary mathematics has involved numerous conceptual innovations, each of which might be seen as paraglorious, and, when taken collectively, might represent a sorites-like inglorious revolution.

So where do Mizrahi and I agree and disagree? He disputes whether there are revolutions exhibiting (TI) in science; I have argued that such revolutions can be found in mathematics.
I take it that we both agree with Crowe that it is a misconception that `the methodology of mathematics is radically different from that of science' \cite[271]{Crowe88}. So we should both like for the story we tell about revolutions to hold for both science and mathematics. Of course, you don't always get what you want---conventional wisdom might suggest that we are both wrong across the board: science and mathematics are methodologically discontinuous in part because science exhibits (inglorious) revolutions but mathematics does not.
Conversely, someone might defend the contrarian stance that Mizrahi and I are both right about revolutions but wrong about the methodological continuity of science and mathematics, because (inglorious) revolutions are confined to mathematics. (This is not as absurd as it may appear---some of the strategies for minimizing the revolutionary aspects of conceptual shifts in science, such as finding common referents between theories, may not work in a field where all the referents are abstract objects.)
However, I believe a more satisfactory resolution is possible.

To see how this might be accomplished, it will help to recast Mizrahi's arguments in my terminology. (TI) may be understood as saying that no revolutions are glorious. So a rebutting defeater against (TI) would be a glorious revolution. Since there are glorious revolutions, Mizrahi concludes that (TI) lacks strong inductive support.
However, if we restrict (TI) to inglorious and paraglorious revolutions, then glorious revolutions no longer count as rebutting defeaters.
%
%Paraglorious and inglorious revolutions should be treated alike since both entail (TI).
%There are plenty of paraglorious revolutions in mathematics---there may even be some inglorious revolutions---but paraglorious revolutions are sufficient to ensure there are many positive instances of (TI).
%There are also many instances of glorious revolutions. Thus Mizrahi's knockdown argument may seem to hold here too: mathematical glorious revolutions would be rebutting defeaters to (TI) with respect to mathematics, since they are revolutions without a significant shift of lexical taxonomy, just as paraglorious or inglorious revolutions would be positive instances.
%But this is not so. Glorious revolutions would be rebutting defeaters to (TI) as a thesis concerning \emph{all} revolutions, but they cannot speak to the claim that all non-glorious revolutions exhibit (TI).
%And this point applies to Mizrahi's argument in science as well as mathematics.
%
%So a resolution: 
Mizrahi does consider, and reject, a related proposal from the `friends of (TI)': retreating to the claim that `some episodes of scientific change exhibit TI, whereas others do not' \cite[372]{Mizrahi14b}.
He rightly objects that such a claim would have no explanatory or predictive value. However, my proposal is more robust: rather than just exclude the anomalous cases, I have offered an independent characterization of subtypes of revolution for which (TI) still holds.% 
\footnote{The contrast is analogous to that Imre Lakatos draws between monster-barring and exception-barring \cite[29]{Lakatos76}.}
Hence (TI) is false as a claim about revolutions in general, as Mizrahi rightly observes. But it is true of two important subtypes: paraglorious and inglorious revolutions.%
\footnote{I am grateful to Moti Mizrahi for the invitation to contribute to this volume and for his insightful comments. I presented an earlier draft at the University of Nevada, Las Vegas: my thanks  for a valuable discussion to the audience, and especially Ian Dove, Maurice Finocchiaro, and James Woodbridge.}

\end{document}